\documentclass[fontsize=11pt,a4paper,DIV13]{scrartcl}

\usepackage[utf8]{inputenc} 
\usepackage{bibitemsep}
\usepackage{cmbright} %PerM: getting error, commenting out!
\usepackage{url}

\usepackage{graphicx}
\usepackage{wrapfig}

\def\NI{\noindent}
\def\ni{\noindent}
\def\sk{\smallskip}
\def\ms{\medskip}
\def\bs{\bigskip}

\def\Proof{\par{\NI\bfseries Proof. }}
\def\qed{\hfill\fbox{\hbox{}}}
\def\qedclaim{\hfill$\triangle$}

\def\term#1{{\em #1}}

\let\it=\itshape
\let\sc=\scshape

\def\SC{\fontfamily{cmr}\fontshape{sc}\selectfont}
\usepackage{amssymb}
\usepackage{amsmath}
\usepackage{amsthm}
\newtheoremstyle{meiner} %name
    {4pt}{3pt}           %space above/below
    {\itshape}           %text font
    {}                   %heading indent
    {\sffamily\bfseries} %heading font
    {}                   %punctuation
    { }                  %space after head
    {}                   %head spec ??
\theoremstyle{meiner}

\newtheorem{theorem}{Theorem}
\newtheorem{lemma}{Lemma}

\newtheorem{proposition}{Proposition}

\def\Case#1.{\rih{Case #1.}}
\def\Phase#1.{\NI{\bfseries Phase #1.}}
\def\Fact#1.{\par\sk{\NI\bfseries Fact #1.}}
\def\Claim#1.{\NI{\bfseries Claim #1.}}

\def\forget#1{}

\def\ITEMMACRO #1 ??? #2 ???{\par\medskip\noindent%
\hangindent=#2em\setbox0\hbox{#1 \kern5pt}%
\ifdim\wd0<\hangindent\setbox0\hbox to\hangindent{\hss#1\quad}\fi%
\box0\ignorespaces}

\def\Item#1{\ITEMMACRO #1 ??? 2.5 ???}

\def\Bitem{\Item{\hss$\bullet$}}

\def\abstract{\noindent\hfil\vbox\bgroup\hsize=.9\hsize
\small\noindent{\bfseries Abstract.}
}
\def\endabstract{\egroup\hfil}

\graphicspath{{Figures/}}
\input{psfig.sty}
\usepackage{subfig}
\usepackage{color}
\def\RR{\hbox{\sffamily I\kern-1ptR}}
\def\NN{\hbox{\sffamily I\kern-1ptN}}
\def\ZZ{\hbox{\sffamily Z\kern-4ptZ}}

\def\dom{}%%  %% _{\mathsf{dom}}}
\usepackage{bbold} %PerM: getting error - commenting out!
\def\onevec{{\mathbb{1}}}
\def\RR{\mathbb{R}}
\def\NN{\mathbb{N}}
\def\ZZ{\mathbb{Z}}
\def\CC{\mathcal{C}}

\def\PsFig#1#2#3{}
\begin{document}
\title{The Complexity of the Partial Order Dimension Problem 
        -- Closing the Gap}
\author{
\parbox{4.6cm}{\center
{\SC Stefan Felsner%%
\footnote{Partially supported by DFG grant FE-340/7-2 and 
              ESF EuroGIGA projects Compose and GraDR.}}\\[3pt]
\normalsize
%%\email{felsner@math.tu-berlin.de}\\
\small
        {Institut f\"ur Mathematik\\
         Technische Universit\"at Berlin\\
         Strasse des 17.~Juni 136\\
         D-10623 Berlin, Germany}}
\and
\parbox{4.4cm}{\center
{\SC Irina Musta\c{t}\u a}\\[3pt]
\small
      {Berlin Mathematical School\\
         Technische Universit\"at Berlin\\
         Strasse des 17.~Juni 136\\
         D-10623 Berlin, Germany}}
\and
\parbox{4.8cm}{\center
{\SC Martin Pergel%
\footnote{
Partially supported by a Czech research grant GA\v{C}R GA14-10799S.}}\\[3pt]
\normalsize
%%\email{perm@kam.mff.cuni.cz}\\
\small
        {Department of Software and Computer Science Education\\
         Charles University\\
         Malostransk\'e n\'am. 25\\ 
         118 00 Praha 1, Czech Republic}}
}%%END AUTHOR

\date{\vskip-14mm}

\maketitle

\subsection*{Abstract}
The dimension of a partial order $P$ is the minimum number of linear
orders whose intersection is $P$.  There are efficient algorithms to
test if a partial order has dimension at most $2$.  In 1982
Yannakakis~\cite{y-cpodp-82} showed that for $k\geq 3$ to test if a
partial order has dimension $\leq k$ is NP-complete. The height of a
partial order $P$ is the maximum size of a chain in $P$.  Yannakakis
also showed that for $k\geq 4$ to test if a partial order of height $2$ has
dimension $\leq k$ is NP-complete. The complexity of deciding whether
an order of height $2$ has dimension $3$ was left open. This question
became one of the best known open problems in dimension theory for
partial orders. We show that the problem is NP-complete.

Technically, we show that the decision problem (3DH2) for dimension is
equivalent to deciding for the existence of bipartite triangle containment
representations (BTCon). This problem then allows a reduction from a class of
planar satisfiability problems (P-3-CON-3-SAT(4)) which is known to be
NP-hard.

\vskip10pt\ni {\bfseries Mathematics Subject
  Classifications (2010)
%  05C10, % Planar graphs / Geometric and topological aspects of graph theory
%  05C62, % Graph representations / geom. and intersection repr.
%  52C15. % Packing and covering in 2D
%  94C15  % Applications of graph theory
% 06A05, % Total order, 
06A07, % Combinatorics of partially ordered sets, 
68Q25, % Analysis of algorithms and problem complexity
05C62, % Graph representations / geom. and intersection repr.
}

\noindent
\textbf{Keywords :} NP-complete, triangle containment, intersection graph. 

% *************************************************************
\section{Introduction}

Let $P = (X, \leq_P)$ be a partial order. A linear order $L = (X,
\leq_L)$ on $X$ is a {\itshape linear extension} of $P$ when $x \leq_P
y$ implies $x \leq_L y$.  A family ${\cal R}$ of linear extensions of
$P$ is a {\itshape realizer} of $P$ if $P = \bigcap_{i} L_i$, i.e., $x
\leq_P y$ if and only if $x \leq_{L} y$ for every $L \in {\cal
  R}$. The {\itshape dimension} of $P$, denoted $\dim(P)$, is the
minimum size of a realizer of $P$.  This notion of dimension for
partial orders was defined by Dushnik and Miller~\cite{dm-pos-41}. The
dimension of $P$ can, alternatively, be defined as the minimum $t$
such that $P$ admits an order preserving embedding into the product
order on $\RR^t$, i.e., elements $x\in X$ have associated vectors
$\hat{x}=(x_1,\ldots,x_t)$ with real entries, such that $x \leq_P y$ if and
only if $x_i \leq y_i$ for all $i$. In the sequel, we denote this
by $\hat{x}\leq\dom\hat{y}$.  The concept of dimension plays a role
which in many instances is analogous to the chromatic number for
graphs.  It has fostered a large variety of research about partial
orders, see e.g.~Trotter's monograph~\cite{t-cpos-92}.

The first edition of Garey and Johnson~\cite{GJ79} listed the decision
problem, whether a partial order has dimension at most $k$, in their
selection of twelve important problems which were not known to be
polynomially solvable or NP-complete. The complexity status was
resolved by Yannakakis~\cite{y-cpodp-82}, who used a reduction from
3-colorability to show that the problem is NP-complete for every fixed
$k\geq 3$. He also showed that the problem remains NP-complete for
every fixed $k\geq 4$ if the partial order is of height~$2$.  The
recognition of partial orders of dimension $\leq 2$ is easy. Efficient
algorithms for the problem have been known since the early 1970's
(e.g.~\cite{epl-pgtg-72}) a linear time algorithm is given
in~\cite{mcs-ltto-97}.  The gap that remained in the complexity
landscape was that for partial orders of height $2$, the complexity of
deciding if the dimension is at most $3$ was unknown. This was noted
by Yannakakis~\cite{y-cpodp-82}, but also listed as Problem~1 in Spinrad's
account on dimension and algorithms~\cite{s-da-94}. It was also
mentioned at several other places, e.g.~\cite{t-gpos:rrnd-96}. In this
paper we prove NP-completeness for this case.

Schnyder's characterization of planar graphs in terms of order
dimension promoted interest in this gap in the complexity landscape.
With a finite graph $G=(V,E)$ we associate the \textit{incidence
  order} $P_{VE}(G)$ of vertices and edgges. The ground set of $P_{VE}(G)$ is $V\cup E$. The
order relation is defined by $x<e$ in $P_{VE}$ if $x\in V$, $e\in E$ and
$x\in e$. The incidence order of a graph is an order of height two.
Schnyder~\cite{s-pgpd-89} proved: A graph $G$ is planar if and
only if the dimension of its incidence order is at most~$3$. This
characterizes a reasonably large class of orders of height 2 and dimension 3.
Efficient planarity tests yield polynomial time recognition algorithms for 
orders in this class. Motivated by Schnyder's
result, Trotter~\cite{t-gpos:rrnd-96} asked about the complexity of
deciding if $\dim(P) \leq 3$ for the class of orders $P$ of height $2$
with the property that every maximal element covers at most $k$
minimal elements.  As a by-product of our proof, we provide an answer:
the problem remains NP-complete even if each element is comparable to
at most $5$ other elements.
\smallskip

In ~\cite{t-pnddtfpos-91} Trotter stated the following interesting problem:

\Item{a)} For fixed $t \geq 4$, is it
NP-complete to determine whether the dimension of the incidence order
of a graph is at most $t$?
\smallskip

\noindent
Dimension seems to be a particularly hard NP-complete problem.  This
is indicated by the fact that we have no heuristics or approximation
algorithms to produce realizers of partial orders that have reasonable
size. A hardness result for approximations was first obtained by Hegde
and Jain~\cite{hj-hapd-07} and recently strengthened by Chalermsook et
al.~\cite{cln-gpr-13}.  They show that unless NP = ZPP there exists no
polynomial algorithm to approximate the dimension of a partial order
with a factor of $O(n^{1-\varepsilon})$ for any $\varepsilon> 0$,
where $n$ is the number of elements of the input order.  The reduction
indeed shows more: 1)~the same hardness result holds for partial
orders of height~$2$.  2)~Approximating the fractional dimension
(c.f.~\cite{bs-fdpo-92,ft-fdpos-94}) of a partial order is also hard.
Another class of partial orders where computing the dimension was
shown to be NP-complete~\cite{kp-cdNfoNPc-89} are the $N$-free partial
orders, i.e., orders whose cover graph has no induced path of length
four.  A motivation for the investigation was that an algorithm for
computing the dimension of $N$-free partial orders would have implied
a constant factor approximation algorithm for general orders.  
As far as approximation is concerned, the best known result is a factor
$O(n\frac{\sqrt{\log\log n}}{\sqrt{\log n}})$ approximation which comes from
an approximation algorithm with the same factor for boxicity, 
see~\cite{abc-approx-12}. 

The following two problems about computational aspects of dimension
are open:

\Item{b)} For fixed $t \geq 3$, is it NP-complete to determine whether
the dimension of an interval order is at most $t$?

\Item{c)} For fixed $w \geq 3$, is it NP-complete to decide whether
the dimension of order of width $w$ is precisely $w$?

% *************************************************************
\subsection{Outline of the paper}
% *************************************************************

In the next section we provide some technical background on 
dimension of orders. In particular we show that $3$-dimensional
orders are characterized by having a containment representation with
homothetic triangles. In Subsection~\ref{ssec:lemmas} we
prove the {\it disjoint paths lemma} and the {\it rotor lemma},
two crucial lemmas for the reduction.

In Section~\ref{sec:reduction} we formally introduce the decision problem
P-3-CON-3-SAT(4), a special planar version of 3-SAT and present the reduction
from P-3-CON-3-SAT(4) to the bipartite triangle containment problem BTCon. We
explain the {\it clause gadget} and the {\it variable gadget} and show how to
build the graph $G_\Phi$ for a given P-3-CON-3-SAT(4) instance $\Phi$. It is
then shown that a triangle containment representation of $G_\Phi$ exists iff there
is a satisfying assignment for $\Phi$. Theorem~\ref{thm:main} summarizes
the result.

In Section~\ref{sec:sandwich} we modify the construction to prove a
`sandwich result', Theorem~\ref{thm:sandwich}. The theorem asserts
hardness of recognition for every class $\CC$ between BTCon and
PUTCon (point unit-triangle containment).  For the proof we introduce
a new variable gadget and construct a graph $H_\Phi$ such that if
$\Phi$ is satisfiable, then $H_\Phi$ has a triangle containment
representation using triangles of only two different sidelengths,
while if $\Phi$ is not satisfiable then $H_\Phi$ has no triangle
containment representation.  As an application we show that
recognition of bipartite graphs that admit a special representation
with triangle intersection (BTInt) is again NP-hard. We then
show that PUTCon and BTInt are in~NP. We conclude
with two open questions and remarks on extensions.

% *************************************************************
\section{Dimension Three and Triangle Containment}
\label{sec:dim+triangles}
% *************************************************************

In the introduction we have mentioned the two definitions of dimension
for partial orders. For completeness, we include a proof
of the equivalence.

\begin{proposition}
A partial order $P=(X,\leq_P)$ has a realizer of size $t$ 
if and only if it admits an order preserving embedding into
$(\RR^t,\leq\dom)$.
\end{proposition}

\Proof
Let $L_1,\ldots,L_t$ be a realizer of $P$. For $x\in X$ and $i\in[t]$ let $x_i$
be the position of $x \in L_i$, i.e., $x_i = |\{ y\;:\; y \leq_{L_i} x\}|$.
The map $x \to (x_1,\ldots,x_t)$ is an order preserving embedding 
into $(\NN^t, \leq\dom)$.  

If an  order preserving embedding of $P$ into $\RR^t$ is given
we consider the projection orders $W_i$ defined 
by $x<_{W_i} y$ if and only if $x_i < y_i$. Elements of $P$ which
share the value of the $i$-th coordinate form antichains
of $W_i$, i.e., being a linear order of antichains $W_i$ is a 
{\itshape weak-order}. If $L_i$ is taken as an arbitrary linear extension
of $W_i$ that is also a linear extension of $P$ this yields a realizer
$L_1,\ldots,L_t$ of $P$.
\qed\bs
 
The definition of dimension via realizers is related to the notion of 
{\itshape reversible sets of incomparable pairs} and {\itshape alternating cycles}
which have been applied throughout the literature on dimension of
partial orders. The definition of dimension via embeddings into product
orders is used less frequently. The most notable application may be
in the context of a generalization of Schnyder's dimension theorem
due to Brightwell and Trotter.

\begin{theorem}[Brightwell-Trotter Theorem~\cite{bt-odcp-93}]
\label{thm:Brightwell-Trotter}
The incidence order $P_{VEF}(G)$ of the vertices, edges and bounded faces of a
3-connected plane graph $G$  has dimension three.
\end{theorem}

Miller~\cite{m-pgmrtmi-02} gave a proof of the theorem by constructing
a rigid orthogonal surface. Based on the connection between Schnyder woods
and rigid orthogonal surfaces, simpler proofs of
Theorem~\ref{thm:Brightwell-Trotter} were obtained in~\cite{f-gepg-03}
and~\cite{fz-swaos-08}. We omit the details here, but
Figure~\ref{fig:graph+rigid-os} shows an example of a rigid orthogonal
surfaces with the corresponding embedding of $P_{VEF}(G)$.

%%%%%%%%%%%%%%%%%%%%%%%%%%%%%%%%%%%%%%%%%%%%%%%%
% in einem figure environment mit caption
   \calc_figscale{27}
    \begin{figure}[htb]
    \centerline{\input{\path/graph+rigid-os.pstex_t}}
    \caption{\label{fig:graph+rigid-os}}
    \end{figure}
    VC
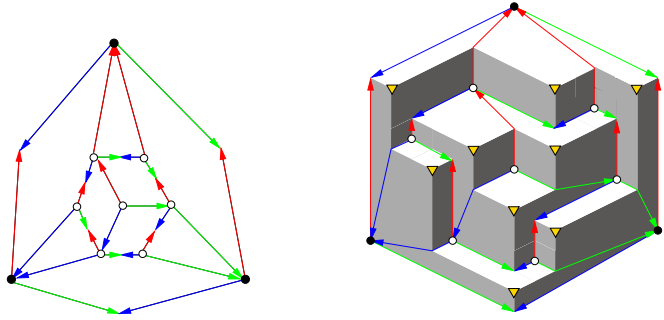
{Left: a 3-connected plane graph $G$ with a Schnyder wood. Right: a
  rigid orthogonal surface corresponding to the Schnyder wood. The
  vertices (circles), faces (triangles) and bend points of edges show
  an order preserving embedding of $P_{VEF}(G)$ into~$\RR^3$.}
%%%%%%%%%%%%%%%%%%%%%%%%%%%%%%%%%%%%%%%%%%%%%%%%
%
The key to our reduction is the characterization of 3-dimensional 
orders which is given in the following proposition.
A {\itshape homothet} of a geometric object $O$ is an object $O'$ that
can be obtained by scaling and translating $O$.
With any family $\cal F$ of sets we associate the
{\itshape containment order} $({\cal F},\subseteq)$.

\begin{proposition}\label{prop:simplex-cont-realizer}
  The dimension of a partial order $P= (X,\leq_P)$ is at most $t$ if and only
  if $P$ is isomorphic to the containment order of a family of homothetic
  simplices in $\RR^{t-1}$.
\end{proposition}

\Proof 
``$\Rightarrow$"
Let $x \to \hat{x}$ be an order preserving embedding of $P$ to
$\RR^t$. With a point $\hat{x}$ in $\RR^t$ associate the (infinite)
orthogonal cone $C(\hat{x}) = \{p \in \RR^t\;:\; p \leq\dom \hat{x}\}$.
Note that $x \leq_P y$ if and only if $C(\hat{x}) \subseteq C(\hat{y})$.
Consider an oriented hyperplane $H$ such that for all $x\in X$ the 
point $\hat{x}$ is in the positive halfspace $H^+$ of $H$, moreover,
$\lambda\onevec\in H^+$ 
and $-\lambda\onevec\in H^-$ for some positive $\lambda$ and the all-ones
vector $\onevec$. For $x\in X$, the intersections of the cone $C(\hat{x})$
with $H$ is a $(t-1)$-dimensional polytope $\Delta(x)$ with $t$ vertices, i.e., 
a simplex. The simplices for different elements are homothetic
and $\Delta(x) \subseteq \Delta(y)$ iff $C(\hat{x}) \subseteq C(\hat{y})$
iff $x \leq_P y$. Hence, $P$ is isomorphic to the containment order
of the family $\{\Delta(x)\;:\; x \in X\}$ of homothetic simplices in
$H \cong \RR^{t-1}$.

\ni ``$\Leftarrow$"
Now suppose that there is a containment order of a family $\cal F$ of
homothetic simplices in $\RR^{t-1}$ that is order isomorphic to $P$.
Let $\Delta(x)$ be the simplex corresponding to $x$ in the order
isomorphism. Apply an affine transformation to get a family ${\cal
  F'} = \{\Delta'(x)\;:\; x \in X\}$ of regular simplices in $H =
\RR^{t-1}$ with the same containment order. Embed $H$ into $\RR^t$ with
normal vector $\onevec$. For each $\Delta'(x)$ there is a unique point
$\hat{x}\in \RR^t$ such that $C(\hat{x}) \cap H = \Delta'(x)$.  Since
the containment orders of $\{C(\hat{x})\;:\; x \in X\}$ and
$\{\Delta'(x)\;:\; x \in X\}$ are isomorphic we identify $x \to
\hat{x}$ as an order preserving embedding of $P$ into $(\RR^t,\leq\dom)$.
\qed\ms

An instance of this proposition, that has been mentioned frequently in the 
literature, is the equivalence between 2-dimensional orders and containment
orders of intervals. 

% *************************************************************
\subsection{Lemmas for triangle containment }\label{ssec:lemmas}
% *************************************************************

From now on we will restrict the attention to partial orders of
height~2.  Note that these orders have a bipartite comparability
graph. Conversely, any bipartite graph with black and white vertices
can be seen as a height~2 order: define $u < w$ whenever $u$ is white,
$w$ is black and $(u,w)$ is an edge. Hence, partial orders of height~2
and bipartite graphs are essentially the same.
  
Given a triangle containment representation of a bipartite graph $G=(V,E)$,
let $B(V)$ be the set of barycenters of the triangles (it can be
assumed that all barycenters are different). Define the
$\beta$-graph $\beta(G)$ as the straight line drawing of $G$ with vertices
placed at their corresponding points of $B(V)$. The following lemma allows
some control on the crossings of edges in $\beta(G)$.

Two edges are called \term{strongly independent} if they share no
vertex (i.e., they are independent) and they are the only edges
induced on their four vertices. The lemma is closely related to the
easy direction of Schnyder's theorem \cite[Theorem 4.1]{s-pgpd-89}.
\goodbreak

\begin{lemma}[disjoint paths lemma]\label{lem:non-cr}
In $\beta(G)$ there is no crossing between two strongly independent edges.
\end{lemma}

\Proof Map the triangle representation via an affine transformation
to a plane $H$ with normal $\onevec$ in $\RR^3$, such that the
triangles on $H$ are equilateral. As in the proof of
Proposition~\ref{prop:simplex-cont-realizer}, we obtain an embedding $v
\to \hat{v}$ of the vertices of $V$ to $\RR^3$ such that the
intersection of the cone $C(v)$ with $H$ is the triangle $T_v$
associated with $v$. The containment relation of triangles $T_v$ on
$H$ and the comparability of points $\hat{v}$ in $\leq\dom$ 
both realize $G$.

\begin{figure}[htb]
    \centerline{\input{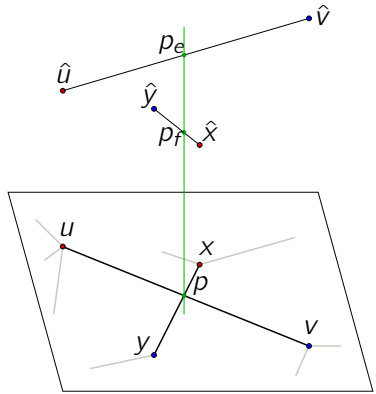}}
    \caption{Illustrating the proof of
  Lemma~\ref{lem:non-cr}.  A crossing pair of strongly
  independent edges in the $\beta$-graph and the implied comparability
  $\hat{x}<\dom \hat{v}$.\label{fig:indep-paths}}
    \end{figure}

Let $e=(u,v)$ and $f=(x,y)$ be strongly independent edges in $G$ such
that $T_u \subset T_v$ and $T_x \subset T_y$. Now suppose that in
$\beta(G)$ the two edges intersect in a point $p$. The line $\ell(p) =
p+\lambda \onevec$ shares a point $p_e$ with the segments $\hat{e} =
[\hat{u},\hat{v}]$ and a point $p_f$ with the segment $\hat{f} =
[\hat{x},\hat{y}]$. Without loss of generality we may assume that
$p_f$ separates $p_e$ and $p$ on $\ell(p)$. Now we look at the
apices of the cones. From triangle containment we obtain
$\hat{u}<\dom p_e <\dom\hat{v}$ and $\hat{x}<\dom p_f<\dom\hat{y}$
and along $\ell(p)$ we get $p<\dom p_f<\dom p_e$ (see Fig.~\ref{fig:indep-paths}).
Combining this and using transitivity we get $\hat{x}<\dom \hat{v}$. 
This implies $T_x \subset T_v$ and
$(x,v)\in E$, contradiction.
\qed
\medskip
%%%%%%%%%%%%%%%%%%%%%%%%%%%%%%%%%%%%%%%%%%%%%%%%
% in einem figure environment mit caption
   \calc_figscale{12}
    
%%%%%%%%%%%%%%%%%%%%%%%%%%%%%%%%%%%%%%%%%%%%%%%%

\begin{figure}[htb]
    \centerline{\input{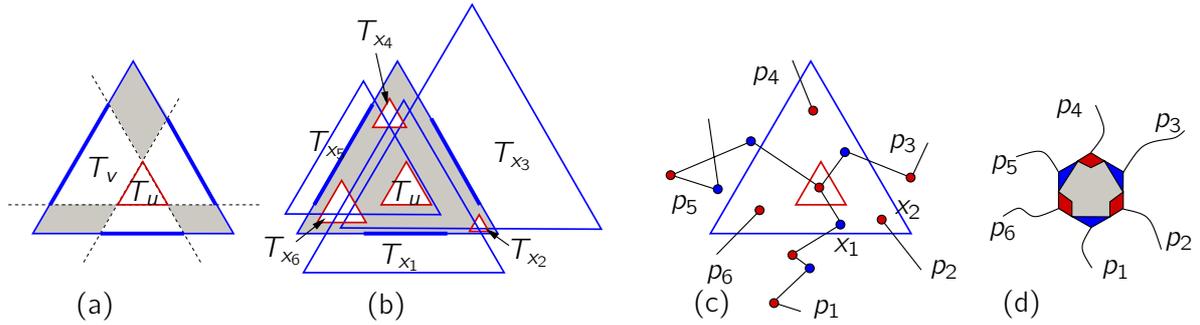}}
    \caption{Part (a) shows the shadow intervals induced by $T_u$ on
  the sides of $T_v$ and in gray the tip regions. Part (b) shows a triangle
  containment representation of a rotor. Part (c) illustrates how the paths $p_1$,
  $p_3$ and $p_5$ together with $u$ in the $\beta$-graph partition the
  interior of $T_v$ in three regions, where the paths $p_2$, $p_4$ and $p_6$
  start. Part (d) shows a schematized picture of an alternating rotor.\label{fig:rotor}}
    \end{figure}

For the reduction we will construct a bipartite graph with an embedding in
the plane which has only few crossings.  Most of these crossings
between edges occur locally in subgraphs that are named
\term{rotor}. A rotor has a center, which is an adjacent pair $u,v$ of
vertices. Additionally there are some non-crossing paths $p_i$.
 Each path $p_i$ is connected to the center either with an
edge $(x_i,u)$ or with an edge $(x_i,v)$ where $x_i$ is an end-vertex
of the path. The interesting case of a rotor is an \term{alternating
  rotor}. In this case it is possible to add a simple closed curve
$\gamma$ to the picture such that (1) $\gamma$ contains the center
and (2) there is a collection of six paths intersecting
the curve $\gamma$ cyclically so that paths leading to $u$ and paths leading
to $v$ alternate.  Figure~\ref{fig:rotor}(b) shows a triangle
containment representation of an alternating rotor.

We will show that in the $\beta$-graph of a triangle containment
representation of an alternating rotor the six paths $p_1,\ldots,p_6$
which do the alternation have to connect to the center in a very
specific way. We assume that $p_1,p_3,p_5$ connect to $u$ and
$p_2,p_4,p_6$ connect to $v$.

Let $u,v$ be the center of a rotor with $T_u \subset T_v$. Define the
\term{shadow intervals} as the intervals between the intersections of the
sides of $T_v$ with the supporting lines of $T_u$, see
Figure~\ref{fig:rotor}(a). As \term{tip regions} of the triangle $T_v$, we
define the three parallelograms determined by the sides of $T_v$ and the
support lines of $T_u$.  We have the following
\medskip

\ni
\Claim 1. In a representation of an alternating rotor
\Bitem 
each of  $T_{x_1},T_{x_3},T_{x_5}$ contains exactly one of the shadow intervals
and each of the shadow intervals is contained in one of the triangles, and
\Bitem
each of $T_{x_2},T_{x_4},T_{x_6}$ intersects exactly one of the three tip
regions of $T_v$ and each of the tip regions is intersected by one of the
triangles.  
\medskip

\Proof From Lemma~\ref{lem:non-cr} we deduce that in the $\beta$-graph
the paths $p_i$ are pairwise non-intersecting, moreover the edges
$(u,x_i)$, with odd $i$, are disjoint from paths $p_j$ with
even~$j$. Hence, the three paths emanating from $u$ define three
regions and each of these regions contains one of the remaining three
paths, see Figure~\ref{fig:rotor}(c). Now each of
$T_{x_1},T_{x_3},T_{x_5}$ contains $T_u$ but it is not contained in
$T_v$, hence, it has to contain a shadow interval. If two of them, say
$T_{x_1}$ and $T_{x_3}$, contain the same shadow interval then there
is no space for $T_{x_2}$ in $T_v$ such that in the $\beta$-graph
$p_2$ is between $p_1$ and $p_3$.  Since $T_{x_j}$, with even~$j$, is
not contained in any $T_{x_i}$, with odd~$i$, it has to intersect one
of the tip regions of $T_v$. Again the alternation of the paths
connecting to $u$ and $v$ implies the one-to-one correspondence
between $T_{x_2},T_{x_4},T_{x_6}$ and the tip regions intersected by
them.  \qedclaim

The essence of the claim is that at an alternating rotor there are six
{\em ports} where paths can attach, these ports alternatingly belong to
the center vertices $u$ and $v$. This property is captured by the
schematic picture of an alternating rotor given in
Figure~\ref{fig:rotor}(d).

%%%%%%%%%%%%%%%%%%%%%%%%%%%%%%%%%%%%%%%%%%%%%%%%
% in einem figure environment mit caption
   \calc_figscale{22}
    
%%%%%%%%%%%%%%%%%%%%%%%%%%%%%%%%%%%%%%%%%%%%%%%%

An \term{alternating 8-rotor} is a rotor with an 8 alternation, i.e.,
there are eight paths $p_1,\ldots,p_8$ intersecting a simple closed
curve $\gamma$ that contains the center cyclically in the order of
indices such that paths leading to $u$ and paths leading to $v$
alternate.

\begin{lemma}\label{lem:8-rotor}
  There is no triangle containment representation of an alternating
  8-rotor such that in the $\beta$-graph of the representation the
  eight paths $p_1,\ldots,p_8$ intersect a simple closed curve around
  the center of the rotor in the order of the indices.
\end{lemma}

\Proof Assume that there is a representation. By considering the
representation of the two alternating rotors with paths $p_1,\ldots,p_6$ and
$p_3,\ldots,p_8$, respectively, we conclude 
from Claim~1 that $T_{x_1}$ and $T_{x_7}$
contain the same shadow interval on $T_v$. This contradicts the 
representation of the alternating rotor with paths $p_1,p_2,p_3,p_4,p_7,p_8$.
\qed

% *************************************************************
\section{The Reduction}\label{sec:reduction}
% *************************************************************

% *************************************************************
\subsection{Decision problems}
% *************************************************************

\ni{\SC Dimension 3 for Height 2 Orders (3DH2)}\\[1mm]
\textbf{Instance:} A partial order $P=(X,<)$ of height 2.\\[1mm]
\textbf{Question:} Is the dimension of $P$ at most 3?
\bigskip

\ni{\SC Bipartite Triangle Containment Representation (BTCon)}\\[1mm]
\textbf{Instance:} A bipartite graph $G=(V,E)$.\\[1mm]
\textbf{Question:} Does $G$ admit a containment representation with
a family of homothetic triangles?
\bigskip

From the $t=3$ case of Proposition~\ref{prop:simplex-cont-realizer} it follows
that the decision problems 3DH2 and BTCon are polynomially equivalent.  In
this section we design a reduction from P-3-CON-3-SAT(4) to show that
BTCon and, hence, also 3DH2 is NP-complete.

In order to define this special version of 3-SAT, we have to recall
the notion of the \term{incidence graph} of a SAT instance
$\Phi$. This is a bipartite graph whose vertices are in correspondence
to the clauses on one side of the partition, and to the variables of
$\Phi$ on the other side. Edges of the incidence graph correspond to
membership of a variable in a clause. Lichtenstein~\cite{l-pfu-82}
showed that the satisfiability problem for CNF-formulae with a planar
incidence graph is NP-complete. Based on this, more restricted
Planar-SAT variants have been shown to be hard. The NP-hardness of
P-3-CON-3-SAT(4) was shown by Kratochv\'\i l~\cite{k-spspc-94}. One of
the first applications was to show the hardness of recognizing
grid intersection graphs.
\bigskip

\ni{\SC P-3-Con-3-SAT(4)}\\[1mm]
\textbf{Instance:} A SAT instance  $\Phi$ with the additional properties:
\Bitem The incidence graph of  $\Phi$ is 3-connected and planar.
\Bitem Each clause consists of 3 literals.
\Bitem Each variable contributes to at most 4 clauses.
\vskip2.5mm
\ni
\textbf{Question:} Does $\Phi$ admit a satisfying truth assignment?
\bigskip

The advantage of working with 3-connected planar graphs, instead of just
planar graphs, is that in the 3-connected case the planar embedding 
is unique up to the choice of the outer face.
Our reduction is inspired by the reduction from~\cite{k-spspc-94}, and even
more so by the recent NP-completeness proof for the recognition of unit grid
intersection graphs with arbitrary girth~\cite{mp-ugig:rp-13}.

% *************************************************************
\subsection{The idea for the reduction}
% *************************************************************

Let $\Phi$ be an instance of P-3-CON-3-SAT(4) with $\Phi$ we assume a
fixed embedding of the incidence graph $I_\Phi$ in the plane. Such an
embedding can be part of the input, otherwise it can be computed
efficiently.  Our aim is to construct a bipartite graph $G_\Phi$ such
that $G_\Phi$ has a triangle containment representation if and only if
$\Phi$ is satisfiable.  The construction of $G_\Phi$ is done by
replacing the constituents of $I_\Phi$ by appropriate gadgets. The
connections between the gadgets are given by pairs of paths that
correspond to the edges of the 3-connected planar graph $I_\Phi$.

The $\beta$-graph of a triangle containment representation of $G_\Phi$
is a drawing of $G_\Phi$. Because of Lemma~\ref{lem:non-cr} all the
crossings in this drawing are local.  In particular there is no
interference (crossings) between different gadgets, no interference
between a gadget and a connecting paths (except when the corresponding
vertex and edge of $I_\Phi$ are incident), and no interference between
two of the connecting paths.  Except for crossings along a path, see
e.g.~the path $p_5$ in Figure~\ref{fig:rotor}(c), we can find
crossings of the drawing only within the gadgets. In particular the
global structure of the drawing, i.e., the cyclic order in which pairs
of paths leave a clause, is as prescribed by $I_\Phi$.

In $G_\Phi$ every edge of the incidence graph
$I_\Phi$ is replaced by a pair of paths. These
two paths join a variable gadget and a clause gadget.  Such a pair of
paths is called an \term{incidence strand}. At their clause ends, the
two paths of an incidence strand are two of the $p_i$ paths of a
rotor. At this rotor the two paths are incident to different center
vertices, hence, they are distinguishable and we may think of one of
them as the green path and the other as the yellow path.  Assuming a
triangle containment representation we look at the two paths of an
incidence strand in the $\beta$-graph. Since they do not cross each
other they define a strip. Looking along this strip in the direction
from the variable gadget to the clause gadget we either see the green
path on the left and the yellow path on the right boundary or the
other way round. This yields an `orientation' of the incidence
strand. The orientation is used to transmit the truth assignment from
the variable to the clause.

The notion of oriented strands is crucial for the design of 
the clause gadgets and variable gadgets described 
in the following two subsections.

% *************************************************************
\subsection{The clause gadget}\label{ssec:clause}
\begin{figure}[htb]
    \centerline{\input{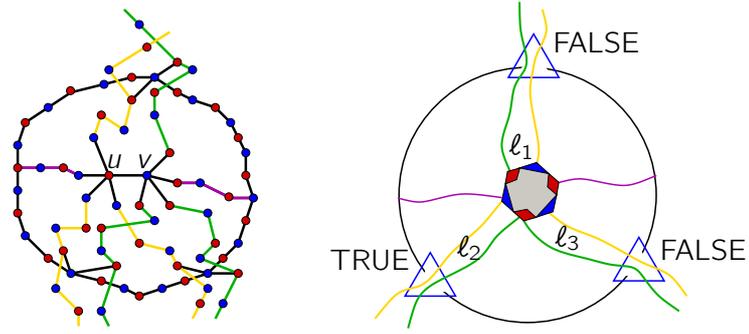}}
    \caption{Left: The graph of a clause gadget.
Right: A schematic view for the case where literals evaluate to $(F,T,F)$.\label{fig:clause-g2}}
    \end{figure}

% *************************************************************
The clause gadget consists of a rotor surrounded by a cycle and two
paths that fix the rotor in the interior of the cycle (magenta). The
paths of the three incidence strands which lead to the clause are also
connected to the rotor. Three vertices of the cycle have two extra
edges connecting to the two paths of an incidence strand. This
construction, we call it \term{gate}, enables the incidence strands to
enter the interior of the cycle. Figure~\ref{fig:clause-g2} shows the
clause gadget as a graph and in a more schematic view.  The lengths of
paths shown in the left figure are chosen at will. For a precise
description of the gadget the lenght of the paths and of the enclosing
circle have to be chosen so that a triangle containment representation
is possible in all cases of a satisfied clause, i.e., in all cases
where the paths are non-crossing. It is quite obvious that if these
lenghts are given by sufficiently large constants, then the triangle
containment representation exists. A reasonable choice could be 314
for the enclosing circle and 50 for each of the eight paths from the
rotor to the circle.  Another detail that has to be verified is that
the gates can be realized via triangle containment.  A possible
solution to this easy exercise is shown in the right part of
Figure~\ref{fig:corner2}.

\begin{figure}[htb]
    \centerline{\input{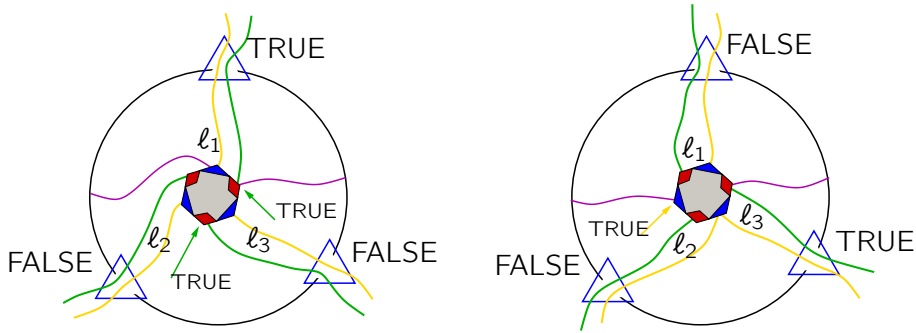}}
    \caption{The left part of the figure shows that 
if the literals evaluate to $(T,F,F)$, then the clause gadget 
can be drawn with noncrossing paths. In fact, by attaching the 
green paths of $\ell_2$ and $\ell_3$ where the arrows show, we
get non-crossing representations for $(T,*,*)$. The right
part shows the same for $(F,*,T)$.\label{fig:clause-s2}}
    \end{figure}

In Figures~\ref{fig:clause-g2} and~\ref{fig:clause-s2} the two paths
connecting the rotor to the cycle are shown in magenta.  The gates are shown
as blue elements (triangles) adjacent to a red vertex of each of the two
paths of the incidence strand associated to the gate. The coloring of the
paths of the incidence strand has been chosen such that at the rotor the
yellow path ends at the red central element $u$ and the green path at the blue
central element $v$.

%%%%%%%%%%%%%%%%%%%%%%%%%%%%%%%%%%%%%%%%%%%%%%%%%%%%%%
% in einem figure environment mit caption
   \calc_figscale{20}
    
%%%%%%%%%%%%%%%%%%%%%%%%%%%%%%%%%%%%%%%%%%%%%%%%%%%%%%

The orientation of the incidence strand corresponding to a literal transmits
TRUE if one of the two paths of the strand can share a port with one of the
magenta paths, it transmits FALSE if the two paths together with an adjacent
magenta path have to use three different ports of the rotor.  To exemplify this
rule: the literal $\ell_1$ is true if the corresponding incidence strand has
yellow left of green when seen from the rotor, whereas the other two literals
are true when on their strand green is left of yellow.  Note that due to this
asymmetry we have to represent a clause as an ordered 3-tuple of literals.

In the right part of Figure~\ref{fig:clause-g2}, we see that 
if the literals evaluate to $(F,T,F)$, then the clause gadget 
can be drawn with noncrossing paths.

Figure~\ref{fig:clause-s2}
shows that in all other cases where the clause evaluates to true,
the schematic clause gadget can as well be drawn without crossing paths,
thus warranting a corresponding triangle containment representation.

%%%%%%%%%%%%%%%%%%%%%%%%%%%%%%%%%%%%%%%%%%%%%%%
% in einem figure environment mit caption
   \calc_figscale{20}
    
%%%%%%%%%%%%%%%%%%%%%%%%%%%%%%%%%%%%%%%%%%%%%%%%

Finally, we have to show that in the $(F,F,F)$ case there is no triangle
containment representation for the clause gadget.  In this case the six paths
of the three incidence strands together with the two (magenta) connecting
paths of the rotor form an alternating 8-rotor, i.e, in the cyclic order the
eight paths have to be connected to the center vertices $u$ and $v$
alternatingly. From Lemma~\ref{lem:8-rotor} we know that such a configuration
has no corresponding triangle containment representation.
We summarize the result in a proposition.

\begin{proposition}\label{prop:clause}
  The clause gadget has a triangle
  containment representation if and only if the two paths of the
  incidence strand of at least one literal are in the orientation
  representing TRUE.
\end{proposition}

% *************************************************************
\subsection{The variable gadget}\label{ssec:var-gad}
% *************************************************************

The variable gadget corresponding to a variable $x$ depends on
the number of occurrences of the variable in clauses. 

We begin with the simplest case. This is when the variable $x$ has only two
occurrences. The basic construction is extremely easy in this case.
The variable $x$ corresponds to one of the incidence strands at each
clause that contain an occurrence of $x$. Let $C_1$ and $C_2$ be the
two clauses containing an occurrence of $x$. Connect the four paths of
the two incidence strands emanating from clauses $C_1$ and $C_2$ to a
single pair of paths. However, there are two options, either we
pairwise connect paths of the same color or we connect paths of
different colors. Which of the two options is the right one to
synchronize the choice of truth values for the literals depends on the
position of the literal containing $x$ in the ordered clauses $C_1$
and $C_2$ and it depends on whether the occurrence is negated or not.
The right pairing for the connection can be determined through an easy
`calculation'.  The result for the case of two positive or two
negative literals is shown in the left table below, where P
corresponds to a color preserving combination of the paths and M for a
color mixing combination.  The right table shows how to connect the
paths if exactly one literal is negated. In this case with respect to
the first table the entries P and M are simply exchanged.
 
$$
\begin{tabular}{l|*{3}{c}}
 & $\ell_1$ & $\ell_2$ & $\ell_3$\\
\hline
$\ell_1$  & M & P & P  \\
$\ell_2$  & P & M & M  \\
$\ell_3$  & P & M & M  \\
\end{tabular}\quad\hskip22mm\quad
\begin{tabular}{l|*{3}{c}}
 & $\ell_1$ & $\ell_2$ & $\ell_3$\\
\hline
$\ell_1$  & P & M & M  \\
$\ell_2$  & M & P & P  \\
$\ell_3$  & M & P & P  \\
\end{tabular}$$

We make an example of how to determine one of the entries of the
tables:
Consider the entry $(\ell_1,\ell_2)$
of the second table.
The variable $x$ is the first literal of a clause $C$ and $\neg x$ is the second
literal of a clause $C'$. For $C$ the literal $\ell_1$ is $T$ if at the rotor the
port of the green path of $\ell_1$ is the clockwise successor of the port of
the yellow path. If $x$ is $T$, then the literal $\ell_2$ of $C'$ has to be
$F$, i.e., at the rotor of $C'$ the port of the green path of $\ell_2$ has to
be a clockwise successor of the port of the yellow path.  Hence, if $x$ is
$T$, then reading the combined strand in clockwise order starting at the rotor
of $C$ is yellow$(\ell_1,C)$, green$(\ell_2,C')$, yellow$(\ell_2,C')$,
and green$(\ell_1,C)$. This shows that in this case the paths have to be
combined so that the colors are mixed. 
 
\begin{figure}[htb]
    \centerline{\input{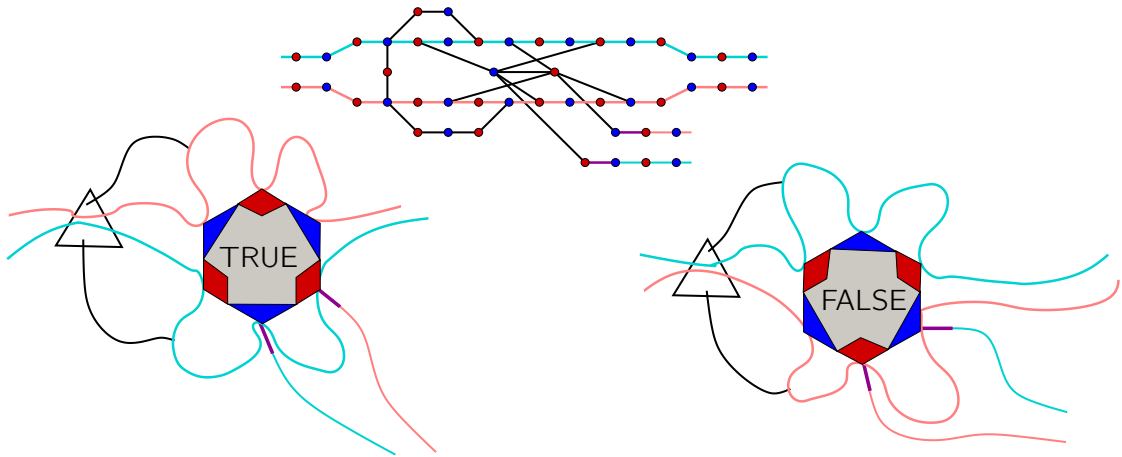}}
    \caption{The variable gadget for a variable with three
  occurrences. The upper picture shows the variable gadget as a bipartite
  graph. The schematic drawings below show the two possible states of the
  variable gadget which correspond to the two possible values of the variable.\label{fig:variable2}}
    \end{figure}

Now assume that the variable $x$ has three occurrences. For two of the
occurrences we can combine the paths of the incidence strands as in
the previous case. To transmit the truth information to a third
incidence strand we use a gadget involving an alternating rotor, the
construction is shown in Figure~\ref{fig:variable2}. The colors are
changed to avoid confusion with between the variable-ends and the
clause-ends of strands.

%%%%%%%%%%%%%%%%%%%%%%%%%%%%%%%%%%%%%%%%%%%%%%%
% in einem figure environment mit caption
   \calc_figscale{20}
    
%%%%%%%%%%%%%%%%%%%%%%%%%%%%%%%%%%%%%%%%%%%%%%%%

In the Figure the combined incidence strands of two occurrences are
shown in cyan and pink. The adjacencies of vertices on these paths to
the center vertices of the rotor between them make an alternating
rotor.  The duty of the black path is to shield two of the ports of
the rotor from outer access. The gate allows access for one of the
incidence strands. There remain two ports of the rotor where
the paths of the third occurrence (magenta ends) can connect to the
center of the rotor. Switching the orientation of the first two
occurrences also makes the strand of the third occurrence switch
orientation, i.e., the truth values transmitted by the three strands
are synchronized. Each connection between a strand of the variable
gadget and the respective strand of a clause gadget has again to be
adjusted so that the correct truth value is received at the clause.
Such an adjustment consists in deciding whether the pairing is
(cyan-green,pink-yellow) or (cyan-yellow,pink-green). Which of the
pairings has to be chosen depends on (1) which strand of the variable
gadget is in question, (2) the position of the literal in the clause,
and (3) whether there is a negation involved. We refrain from listing
the possible cases.

Finally, if the variable $x$ has four occurrences, then we take two
of the gadgets used for variables with three occurrences and 
connect two paths of a strand of each of them.
This synchronizes the truth values transmitted by the remaining
four strands. 

% *************************************************************
\subsection{Wrap-up}
% *************************************************************

We have described the gadgets needed to construct the instance
$G_\Phi$ for the decision problem BTCon based on an instance $\Phi$
for P-3-Con-3-SAT(4). The discussion of the clause gadget lead to
the insight that a clause gadget has a triangle containment representation 
if and only if at least one of its literals evaluates to true 
(Proposition~\ref{prop:clause}).

\begin{figure}[htb]
    \centerline{\input{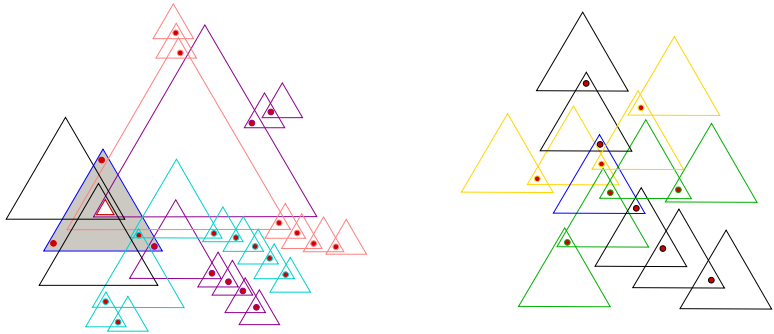}}
    \caption{Left: The rotor of a variable gadget with a
  pink path bypassing a magenta path at a blue port and a cyan path
  bypassing a magenta path at an adjacent red port.\\ 
  Right: A realization of a gate.\label{fig:corner2}}
    \end{figure}

Regarding the variable gadget, we claim that a triangle containment
representation exists whenever all the involved paths are given by a
large enough constant. To see this we need the representation of an
alternating rotor, and we have to verify that it is possible to
realize the ports where a pink or cyan path is bypassing a magenta
path (colors as in Figure~\ref{fig:variable2}). The realizability of
the configuration is shown in Figure~\ref{fig:corner2}. Regarding the
length of the path we propose to use 20 for each segment of a path
between two vertices of degree three. (The value 20 is rather 
arbitrary, we chose it because it allows a realization of the variable
gadged which is clearly arranged.)
\vskip-6mm \vbox{}

%%%%%%%%%%%%%%%%%%%%%%%%%%%%%%%%%%%%%%%%%%%%%%%
% in einem figure environment mit caption
   \calc_figscale{13}
    
%%%%%%%%%%%%%%%%%%%%%%%%%%%%%%%%%%%%%%%%%%%%%%%%

The graph $G_\Phi$ is a subdivision of some core graph
$G_\Phi^\bullet$ with minimum degree 3. To actually build the graph
$G_\Phi$, the length of the paths replacing edges of $G_\Phi^\bullet$
have to be specified. For the clause gadget and for variable gadgets
we identify subdivisions so that they are representable in all cases
except the unsatisfied clause. The corresponding path lengths are
copied to all edges of $G_\Phi^\bullet$ that are internal to a
gadget. For the incidence strands we have initial contributions to the
path length that come from the gadgets. For the final adjustment of
the path lengths of incidence strands, we consider a fixed straight
line drawing of $I_\Phi$ on a small integer grid. The length of the
path of an incidence strand can be taken as a multiple of the edge
length of the corresponding edge in the drawing of $I_\Phi$.  The
proportionality factor should be large enough to ensure that the paths
can be represented with a sequence of small triangles so that unwanted
interference between different paths can be avoided.  If the enclosing
circle of the clause gadget has length 314 as suggested in
Subsection~\ref{ssec:clause}, then we may take the proportionality
factor to be 500. Again the choice of 500 is rather arbitrary, but we
think of a representation where the small triangles are points or
almost points and the large triangles have sidelength close to
one. With this assumption the clause and variable gadgets can live in
disks of radius 25, hence they are small compared to the distance that
can be spanned with a path, resp. strand, connecting two of the
gadgets.

In summary: Given an instance $\Phi$
for P-3-Con-3-SAT(4) we can construct a bipartite graph $G_\Phi$
such that $G_\Phi$ admits a triangle containment representation 
if and only if $\Phi$ is satisfiable.

Since P-3-Con-3-SAT(4) is NP-complete~\cite{k-spspc-94}, this shows that 
BTCon is NP-hard. To show membership of BTCon in the class NP we
recall that BTCon is equivalent to 3DH2. Membership of 3DH2 in the class NP
is quite immediate, a realizer can serve as certificate.
The certificate can be checked efficiently.

The theorem below claims that BTCon remains NP-complete for bipartite input
graphs with maximum degree at most 5 and arbitrarily large constant
girth. Regarding the degree, it suffices to check that all the vertices
involved in clause and variable gadgets have degree at most 5 (see
Figure~\ref{fig:clause-g2} and~\ref{fig:variable2}).  Indeed, only center
vertices of rotors are of degree~5 and degree~4 only occurs at gates in clause
gadgets. Regarding the girth, we observe that every cycle of $G_\Phi^\bullet$
contains an edge that is subject to subdivisions. Hence the girth can be made
arbitrarily large.

\begin{theorem}\label{thm:main}
The recognition problems BTCon and 3DH2 are NP-complete.
NP-hardness of the recognition is preserved for input graphs of
arbitrarily large (constant) girth and maximum degree~5.
\end{theorem}

% *************************************************************
\section{A Sandwich Theorem}\label{sec:sandwich}
% *************************************************************

The proof of Theorem~\ref{thm:main} in the previous section
was based on 3 properties of triangle containment representations.

\Item{1)} Validity of Lemma~\ref{lem:non-cr} implying that independent
     paths are non-crossing in any representation.

\Item{2)} The restricted representability of an alternating rotor 
      and the non-representability of an alternating 8-rotors 
      (Lemma~\ref{lem:8-rotor}).

\Item{3)} If $G_\Phi$ admits a schematized  drawing without crossings
     and the path lengths are assigned properly, then 
     $G_\Phi$ admits a triangle containment representation.
\medskip

\ni These three properties also hold for restricted classes of
bipartite triangle containment representations. For example, we may
require that the representing triangles are pairwise homothetic and of
only two sizes, i.e., form two equivalence classes under translation.
This is equivalent to the condition that the triangles representing
one color class of the bipartite graph are points and
triangles representing the other class are translates of a fixed
triangle. This fixed triangle can be assumed to be equilateral 
with unit sidelength.
\bigbreak

\ni{\SC Point Unit-Triangle Containment (PUTCon)}\\[1mm]
\textbf{Instance:} A bipartite graph $G=(V,E)$.\\[1mm]
\textbf{Question:} Does $G$ admit a triangle containment representation
with points and translates of an equilateral unit triangle.
\bigskip

The following theorem states that the decision problem PUTCon
is hard. Even stronger: recognition of every class of graphs that is 
sandwiched between PUTCon and BTCon is also NP-hard. 

\begin{theorem}[Sandwich Theorem]\label{thm:sandwich}
The recognition problems PUTCon and BTCon are NP-complete.
Moreover:
\Bitem 
Recognition of every class $\CC$ of bipartite graphs 
such that $\CC$ contains all yes instances of PUTCon and
$\CC$ is contained in the set of all yes instances of BTCon  
is also NP-hard. 
\Bitem 
NP-hardness of the recognition of $\CC$ is preserved for input graphs of
arbitrarily large (constant) girth and maximum degree at most 7.
\end{theorem}

\Proof
The idea for the proof is to encode a P-3-Con-3-SAT(4) instance $\Phi$
in a graph $H_\Phi$ such that:
\Item{(F)}
If there is no satisfying assignment for $\Phi$, then there is no 
BTCon representation for $H_\Phi$. Since  $\CC \subset $ BTCon  this implies 
that $H_\Phi \not\in \CC$.
\Item{(T)}
If there is a satisfying assignment for $\Phi$, then there is a
PUTCon representation for $H_\Phi$. Since  PUTCon $\subset \CC$ this implies 
that $H_\Phi \in \CC$.
\medskip

\ni
Property (F) is satisfied by the graph $G_\Phi$ from the
previous section. Unfortunately the graph $G_\Phi$ has no PUTCon
representation, the problem is that in the variable gadget 
we need that a path can bypass another path at a port of a rotor,
see Figure~\ref{fig:corner2}(Left). To realize a bypass at each of two
consecutive ports, we need triangles of at least three different sizes.

Clause gadgets with sufficiently long
paths can be realized in the PUTCon model provided that one of the
literals represents TRUE. The realization of the gate with points and
unit triangles is shown in Figure~\ref{fig:corner2}(Right).

To make the approach work, we define $H_\Phi$ on the basis of our
clause gadgets from Subsection~\ref{ssec:clause} but with a new
variable gadget.  Since Property (F) depends on the clause gadget, it
is also valid for~$H_\Phi$. 

The new variable gadget has space to accommodate six occurrences of
the variable in clauses.  Since in our instances variables only have
three or four occurrences, we use a {\em waste gadget} (see below) to
get rid of some of them. The variable gadget superimposes two
alternating rotors on the same two center vertices $u$ and $v$. Let
$r_1,..,r_6$ and $l_1,..,l_6$ be, in cyclical order, the six vertices
of the first rotor and second rotor, respectively, that are connected
to the center. However, for each $i$ one of $r_i$ and $l_i$ is
connected to $u$ and the other to $v$. The paths starting at $r_i,l_i$
form an incidence strand or they form a waste strand. The ordering of
the strands around the variable is in accordance with the plane
embedding of $I_\Phi$.

\begin{figure}[htb]
    \centerline{\input{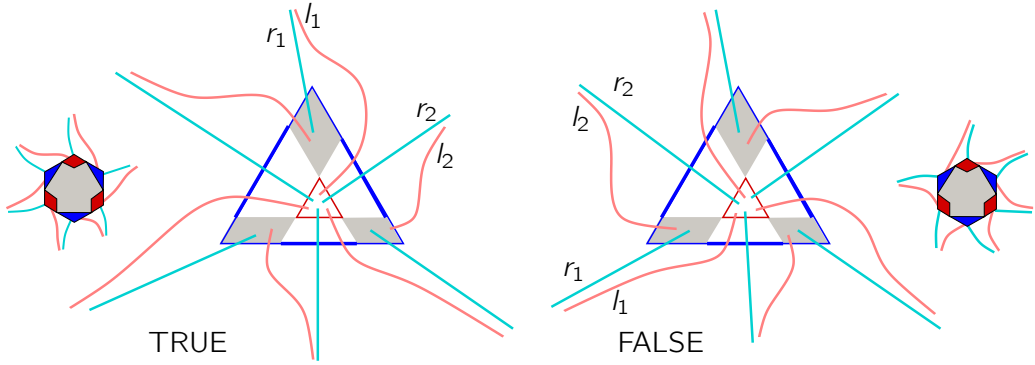}}
    \caption{Two pictures for each of the two possible ways of
  interlacing the rotors of a variable gadget, as defined by the truth value. The
  small pictures use schematized rotors, while the larger pictures show the 
  triangles of the center vertices of the rotors.\label{fig:variable}}
    \end{figure}

For $i\neq j$ there are no crossings between any paths starting at
$l_i$ and $r_j$, hence, the two paths of an incidence strand have to
connect adjacent ports of the rotor. Therefore there are just two
possible ways in which the two rotors can interlace, see
Figure~\ref{fig:variable}.  These two configurations correspond to the
assignment of the values TRUE or FALSE to the variable. As in the case
of the variable gadget from Subsection~\ref{ssec:var-gad}, the strands 
of the variable have to be connected to the strands of the clause
by a proper choice of the pairing 
(cyan-green,pink-yellow) or (cyan-yellow,pink-green).

%%%%%%%%%%%%%%%%%%%%%%%%%%%%%%%%%%%%%%%%%%%%%%%%
% in einem figure environment mit caption
   \calc_figscale{20}
    
%%%%%%%%%%%%%%%%%%%%%%%%%%%%%%%%%%%%%%%%%%%%%%%%

\begin{figure}[htb]
    \centerline{\input{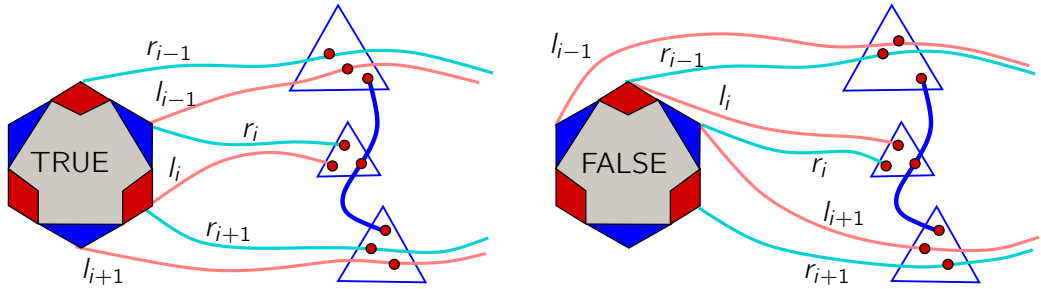}}
    \caption{Two instances of a waste gadget. The hexagon
  belongs to a variable gadget.\label{fig:waste}}
    \end{figure}

The \term{waste gadget} is designed to keep wasted strands in place between
two incidence strands. For this purpose, they are connected to the neighboring
incidence strands as shown in Figure~\ref{fig:waste}.

%%%%%%%%%%%%%%%%%%%%%%%%%%%%%%%%%%%%%%%%%%%%%%%%
% in einem figure environment mit caption
   \calc_figscale{25}
    
%%%%%%%%%%%%%%%%%%%%%%%%%%%%%%%%%%%%%%%%%%%%%%%%

The new variable gadget as well as the waste gadgets have PUTCon
representations, hence, Property (T) is also true for $H_\Phi$.  The claim
regarding the girth follows, because each cycle contains a path that can be
subdivided accordingly, whereas for the maximum degree we only have to note
that the maximum degree in the new variable and waste gadgets is 7.  \qed
 
% *************************************************************
\subsection{An application: bipartite triangle intersection representations}
% *************************************************************

We use Theorem~\ref{thm:sandwich} to show hardness of recognizing
bipartite graphs that admit a special type of triangle intersection
representation.

\ni{\SC Bipartite Triangle Intersection Representations (BTIntR)}\\[1mm]
\textbf{Instance:} A bipartite graph $G=(V,E)$ with bipartition $X,Y$.\\[1mm]
\textbf{Question:} Are there families $\Delta_X = \{T_x\;:\; x \in X\}$ and 
$\Delta_Y = \{T^*_y\;:\; y \in Y\}$ of triangles such that
\Bitem Triangles in each of the families are pairwise homothetic, and
 there is a point reflection transforming $T_x$ into $T^*_y$. 
\Bitem $(x,y) \in E$ if and only if
$T_x \cap T^*_y \neq \emptyset$.
\Bitem Within each of the two families $\Delta_X$ and $\Delta_Y$ there
is no containment.
\bigskip

\ni As a consequence of Proposition~\ref{prop:separated-simplices}, we
obtain that the problem BTIntR is equivalent to the question of
whether a partial order $P$ of height 2 admits an order preserving
embedding $P \to \RR^3$ such that there is a plane $H$ in $\RR^3$ that
separates the minimal and the maximal elements of $P$. In particular,
$P$ has to be 3-dimensional. The partial orders corresponding to
bipartite graphs in PUTCon clearly admit such a `separated' embedding
in $\RR^3$, hence, yes instances of BTIntR form a class of graphs
sandwiched between PUTCon and BTCon. With Theorem~\ref{thm:sandwich}
this implies hardness. In Proposition~\ref{prop:BTIntR in NP} we show
that BTIntR is in NP. Together this proves the following:

\begin{theorem}
The decision problem BTIntR is NP-complete.
\end{theorem}

The proof of Proposition~\ref{prop:separated-simplices} is based on
representations of $t$-dimensional orders using simplices
in~$\RR^{t-1}$. In Proposition~\ref{prop:simplex-cont-realizer} we
have shown that $t$-dimensional orders are containment orders of
homothetic simplices in $\RR^{t-1}$.

Now consider a representation $\hat{X}=\{\hat{x}\;:\; x \in X\}$ of
$P=(X,<_P)$ in $(\RR^t,\leq\dom)$ and let $H$ be a hyperplane with
normal vector $\onevec$ and an element $h_0\in H$. Define the 
half-spaces $H^+=\{y : \langle y-h_0,\onevec\rangle \geq 0 \}$ and 
$H^-=\{y : \langle y-h_0,\onevec\rangle \leq 0 \}$ and two subsets $X^+ = \{x\in
X\;:\;\hat{x}\in H^+\}$ and $X^- = \{x\in X\;:\;\hat{x}\in H^-\}$ of $X$. 

An element $x\in X^+$ is represented on $H$ by the simplex
$\Delta(x) = H \cap C(\hat{x})$, where $C(\hat{x})$ is the cone
$\{p \in \RR^t\;:\; p \leq\dom \hat{x}\}$. The containment
relation on $\{\Delta(x)\;:\; x \in X^+\}$ represents the 
induced order $P[X^+]$.

For elements of $X^-$ we consider the dual
cone $C^*(\hat{y}) = \{p \in \RR^t\;:\; p \geq\dom \hat{y}\}$.  Note
that the containment of dual cones represents the dual $P^*[X^-]$ of
the suborder $P[X^-]$ of $P$ induced by $X^-$, i.e., $C^*(\hat{x})
\subseteq C^*(\hat{y})$ if and only if $y \leq_P x$.  The
intersections of $H$ with the dual cones yields a family homothetic
regular simplices $\{\Delta^*(y)\;:\; y \in X^-\}$ whose containment
order represents $P^*[X^-]$.

If $x\in X^+$ and $y\in X^-$, then
$\Delta(x)$ and $\Delta^*(y)$ are both regular simplices, however,
they are not homothetic, a point reflection is needed to get from one
to the other. Comparabilities between $x\in X^+$ and $y\in X^-$ also
have a nice description in the family $\{\Delta(x)\;:\; x \in
X^+\}\cup \{\Delta^*(y)\;:\; y \in X^-\}$:

\par\sk{\NI\itshape Claim.}  $x\in X^+$ and $y\in
X^-$ are comparable if and only if $\Delta(x) \cap \Delta^*(y)
\neq\emptyset$.

\par\sk\NI
A comparability has to be of the form $y <_P x$. Hence, if
$\overline{xy}$ is the segment with endpoints $\hat{x}$ and~$\hat{y}$,
then $\overline{xy} \subset C(\hat{x}) \cap C^*(\hat{y})$. It follows
that $\overline{xy}\cap H$ is a point in the intersection $\Delta(x)
\cap \Delta^*(y)$.

Conversely, if $p\in \Delta(x) \cap \Delta^*(y)$, then
$\hat{y}\leq\dom p \le\dom \hat{x}$, i.e., $y <_P x$. 
\qedclaim
\ms

\begin{figure}[htb]
    \centerline{\input{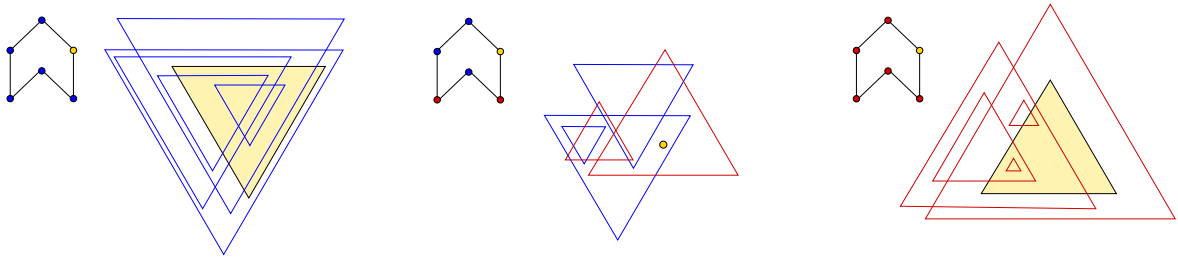}}
    \caption{Three representations of the six element partial order $P$ (the chevron)
with triangles. The representations are based on the same embedding in $\RR^3$ 
with different choices of the hyperplane $H$. Left: $H$ is below all points $\hat{x}$ and we
obtain a containment representation with triangles. Center:
$H$ contains $\hat{x_0}$, where $x_0$ is the element that is distinguished
in all three representations. The simplex of $x_0$ is degenerated to a point
elements of $X^+$ are blue and elements of $X^-$ are red. 
Right: $H$ is above all points $\hat{x}$ and the containment of triangles
represents the dual order $P^*$.\label{fig:chevron-triangles2}}
    \end{figure}

Figure~\ref{fig:chevron-triangles2} shows an order $P$ with
the triangle configurations obtained from different choices of $H$ in
an embedding of $P$ into $\RR^3$.

%%%%%%%%%%%%%%%%%%%%%%%%%%%%%%%%%%%%%%%%%%%%%%%%
% in einem figure environment mit caption
   \calc_figscale{16}
    
%%%%%%%%%%%%%%%%%%%%%%%%%%%%%%%%%%%%%%%%%%%%%%%%

Based on these considerations the following
certificate for dimension at most $t$ is easily obtained.

\begin{proposition}\label{prop:separated-simplices}
  Let $P=(X,\leq_P)$ be a partial order with a partition $X = X^-\cup X^+$
  such that $X^-$ is a downset and $X^+$ is an upset in $P$. If there exist
  families of full-dimensional regular simplices $\{\Delta(x)\;:\; x \in X^+\}$
  and $\{\Delta^*(y)\;:\; y \in X^-\}$ in $\RR^{t-1}$, where the simplices in
  each of the families are pairwise homothetic and there is a point
  reflection transforming $\Delta(x)$ into $\Delta^*(y)$ such that
\Bitem $x \leftrightarrow \Delta(x)$ is an isomorphism between
       $P[X^+]$ and the containment order on $\{\Delta(x)\;:\; x \in X^+\}$
\Bitem $y \leftrightarrow \Delta^*(y)$ is an isomorphism between
       $P^*[X^-]$ and the containment order on $\{\Delta^*(y)\;:\; y \in X^-\}$,
\Bitem $x\in X^+$ and $y\in X^-$ are comparable 
   if and only if $\Delta(x) \cap \Delta^*(y) \neq \emptyset$,
\par\medskip\ni
then $\dim(P) \leq t$.
\end{proposition}

\begin{proposition}\label{prop:BTIntR in NP}
 BTIntR and PUTCon are in NP.
\end{proposition}

\Proof Kratochv\'{\i}l and Matou\v{s}ek~\cite{km-igs-94} showed that
the class of segment intersection graphs where the segments of the
representation are restricted to the elements of a fixed set $K$ of
directions with $|K|=k$ is in NP for every $k$. Our proof is based on
their ideas.  As certificate of a yes instance for the BTIntR problem
we propose a combinatorial description of the intersection pattern of
the triangles. In this description each vertex comes with three lists,
each of them representing the order of intersections along a side of
the triangle. To check whether the certificate corresponds to a
BTInt/PUTCon representation, we can use a linear program to test
whether the combinatorial description can be ``stretched'' to a set of
homothetic triangles with the same intersection pattern. 

The requirement that the intersection of side $s$ with $q$ immediately
precedes its intersection with $r$ is encoded in a linear inequality:
Consider the equations $a_s x + b_s=y$, $a_r x + b_r = y$, and $a_q x
+ b_q = y$ of the supporting lines and note that the $a_i$ are slopes,
i.e., they are constants and can be taken from the fixed 3-element set
$\{-1,0,1\}$, the $b_i$ are the variables. The $x$ coordinate of the
intersection of~$s$ and $r$ is $\frac{b_s-b_r}{a_s-a_r}$. The
condition that this is smaller than the $x$ coordinate of $s\cap q$ is
captured by the strict inequality $(b_s
-b_r)(a_s-a_q)<(b_s-b_q)(a_s-a_r)$. Introducing some constant $c$ for
the min distance, we obtain the inequality $(b_s -b_r)(a_s-a_q) + c
\leq (b_s-b_q)(a_s-a_r)$ in the variables $b_s,b_r,b_q$. In the case
of PUTCon we additionally need that the distance between two corners
of a triangle realizes a prescribed value.  Since corners also are
intersections of segments the same technique applies.

From the system of equations we make linear program by adding the
objective $\max c$.  The decision whether the linear program has a
solution with objective $c > 0$ can be obtained in polynomial
time.\qed

% *************************************************************
\section{Open questions and extensions}\label{sec:concl}
% *************************************************************

%
%In the introduction we have already stated four problems
%about computational aspects of dimension that we consider to
%be interesting. The considerations of this paper suggest 
%further questions.

We have shown that maximum degree at least 5 is enough to make BTCon
hard. From Schnyder's
work~\cite{s-pgpd-89} we know that bipartite graphs
where the degree in one of the color classes is~2 
are yes instances for
BTCon if and only if they are incidence orders of planar graphs.
What about maximum degree 3? 

\Bitem What is the complexity of deciding whether a bipartite graph of
maximum degree 3 admits a BTCon representation?
\smallskip

\ni
We can also restrict the class of inputs to planar bipartite graphs.
It is known that the incidence order of vertices and faces of a
3-connected planar graph is of dimension 4, see~\cite{bt-odcp-93},
moreover there are outerplanar graphs whose incidence order of vertices
and faces is of dimension 4, see~\cite{fn-dopm-07}. Is it hard to decide whether a
planar bipartite graph is of dimension 3? Or in terms
of triangle containments:

\Bitem 
What is the complexity of deciding whether a planar bipartite
graph admits a BTCon representation?
\smallskip

\ni
From Schnyder's proof we know that 
incidence orders of planar graphs have BTInt representations.

\Bitem What is the complexity of deciding whether an incidence orders
of planar graphs (a subdivision of a planar graph) 
admits a PUTCon representation?
\smallskip

\ni
Constructions with rotors with larger or smaller numbers of corners
may yield other NP-hardness proofs. An example is given by the
proof of hardness of recognition of unit grid intersection graphs
(UGIG) in~\cite{mp-ugig:rp-13} where rotors have four corners.
In fact, the paper shows hardness of recognition for all classes that
are sandwiched between UGIG and pseudosegment intersection graphs.

We invite the reader to adapt the rotors and the gadgets to show that
intersection graphs of polygonal regions with $k$ corners are hard to
recognize for every fixed $k$. This has already been shown by
Kratochv\'{\i}l~\cite{k-spspc-94}, however, his proof does not show
hardness for inputs of large girth.

% *************************************************************
%%\break
\def\sc{\SC}\small
\advance\bibitemsep-3pt
\bibliographystyle{my-siam} 
\bibliography{3dim}
% *************************************************************
\end{document}